\documentclass{birkjour}
\usepackage{amsmath, amsfonts, amssymb, graphicx, color}
\usepackage{mathrsfs}
\usepackage{latexsym}

 \newtheorem{thm}{Theorem}[section]
 \newtheorem{cor}[thm]{Corollary}

 \theoremstyle{definition}
 \newtheorem{defn}[thm]{Definition}
 \theoremstyle{remark}
 \newtheorem{rem}[thm]{Remark}

 \numberwithin{equation}{section}
\def\({\left ( }
\def\){\right )}
\def\<{\left < }
\def\>{\right >}

\begin{document}

\title[Characterizations of $\ast$-Ricci-Bourguignon.....]{Characterizations of $\ast$-Ricci-Bourguignon solitons on Kenmotsu manifolds}
\author[S. Roy]{Soumendu Roy}
\address{Department of Mathematics, School of Advanced Sciences, Vellore Institute of Technology, Chennai-600127, India}
\email{soumendu1103mtma@gmail.com, soumendu.roy@vit.ac.in}
\author[K.R]{Karthika Ramasamy}
\address{Department of Mathematics, School of Advanced Sciences, Vellore Institute of Technology, Chennai-600127, India}

\email{karthika.r2023@vitstudent.ac.in}

\author[L.K]{Lavanya Kumar}

\address{Department of Mathematics, School of Advanced Sciences, Vellore Institute of Technology, Chennai-600127, India}

\email{lavanya.k2023a@vitstudent.ac.in}
\author[P. Jana]{Purabi Jana}

\address{Department of Mathematics, School of Advanced Sciences, Vellore Institute of Technology, Chennai-600127, India}

\email{purabi.jana2023@vitstudent.ac.in}



\begin{abstract}
In this paper, we have found some features of $\ast$- Ricci Bourguignon Soliton on Kenmotsu manifold. We estimated the conditions for $\ast$-Ricci Bourguignon on Kenmotsu manifold to be compressing, balancing or enlarging accordingly. We have found some curvature properties of Kenmotsu manifold admitting $\ast$-Ricci Bourguignon Soliton. Additionally, we have featured $\ast$-Ricci Bourguignon on Kenmotsu manifold with torse-forming vector field. Finally, we proved an example of $5$-dimensional Kenmotsu manifold on $\ast$-Ricci Bourguignon Soliton. 
\end{abstract}

\maketitle

\section{Introduction}\label{1}
\hspace{.5cm}
On the basis of  some unpublished work by Lichnerowicz and a paper by Aubin \cite{AU}, Jean-Pierre Bourguignon \cite{PJB} created a new geometric flow in 1981 that he called the Ricci-Bourguignon flow \cite{SD}, which can be demonstrated as, 
\begin{defn}
If a collection of metrics $g(t)$, on an n-dimensional manifold $\mathfrak{M}_{n}$ fulfills the following evolution equation, then $g(t)$ will be expressed as the Ricci- Bourguignon flow [abbreviated as RB flow]\\
\begin{equation}\label{1.1}
 \frac{\partial g}{\partial t}= -2\mathbb{S}-2 \omega r g
\end{equation}
where $\mathbb{S}$ is the Ricci tensor, $\omega$ is a real constant and $r$ is the scalar curvature.
\end{defn}
Shubham Dwivedi \cite{SD} introduced the Ricci-Bourguignon Soliton. A homogeneous solution to RB flow will be termed as a Ricci Bourguignon soliton and it is defined as, If a vector field $z$ on $\mathfrak{M}$(Riemannian Manifold) that assures \cite{DY} that assures\\
\begin{equation}\label{1.2}
 \mathscr{L}_z g +2\mathbb{S}=2\Omega g +2\omega rg
\end{equation}
    where  $\mathscr{L}_z g$ is a Lie derivative of the metric $g$ with respect to the vector field $z$ and $\Omega \in \mathbb{R}$ a constant.\\
    According to $\Omega\gtreqqless0$, it is said to be compressing, balancing or enlarging respectively.\\
    Tachibana \cite{TS} and Hamada \cite{TH} invented the $\ast$-Ricci tensor and it is defined as 
    \begin{equation}\label{1.3}
        \mathbb{S}^{\ast}(\mathbb{Y}_1,\mathbb{Y}_2)=\frac{1}{2}(Tr[\phi,r(\mathbb{Y}_1,\Phi{\mathbb{Y}_{2}})])
    \end{equation}
    for all vector fields $\mathbb{Y}_{1}$,$\mathbb{Y}_{2}$ on $\mathfrak{M}^n$, $\Phi$ is a $(1,1)$- tensor fields.\\
    The above equation $\eqref{1.3}$ is said to develop $\ast$-$\eta$-Einstein manifold, if its satisfies\\
    $$\mathbb{S}^{\ast}(\mathbb{Y}_1,\mathbb{Y}_2)=\lambda{g}(\mathbb{Y}_1,\mathbb{Y}_2)+ \mathcal{C} \eta(\mathbb{Y}_{1})\eta(\mathbb{Y}_{2})$$
    for all vector fields $\mathbb{Y}_{1}$,$\mathbb{Y}_{2}$ and $\lambda$, $\mathcal{C}$ are smooth functions.\\
    In addition, if $\mathcal{C}$ =0, then\\
    $$\mathbb{S}^{\ast}(\mathbb{Y}_1,\mathbb{Y}_2)=\lambda{g}(\mathbb{Y}_1,\mathbb{Y}_2)$$
    for all vector fields $\mathbb{Y}_1$,$\mathbb{Y}_2$ then the manifold\cite{RA} will be reduced as $\ast$- Einstein.\\
 
The Riemannian manifold $(\mathfrak{M}^n,g)$ is said to be a $\ast$- Ricci Bourguignon Soliton (abbreviated as $\ast$-RB Soliton ) \cite{SD1} if there is a vector field $z$ on $M$ its satisfies 
 \begin{equation}\label{1.4}
     \mathscr{L}_z g+ 2 \mathbb{S}^{\ast} = 2 [\Omega + \omega r ^{\ast}] g
 \end{equation}
 here, $\mathscr{L} _z g $ is the Lie derivative of the metric $g$ relative to the vector field $z$, $\Omega$ $\in$ $\mathbb{R}$ is said to be constant , $\mathbb{S}^{\ast}$ will be $\ast$- Ricci tensor and $r^{\ast}=Tr(\mathbb{S}^{\ast})$ and finally it is known to be $\ast$- scalar curvature.As well as, it is said to be compressing, balancing or enlarging regarding 
 $\Omega\gtreqqless 0$ accordingly.  If in $\eqref{1.4}$ $\Omega$ is smooth function but not as a constant, then $\ast$-Ricci-Bourguignon Soliton \cite{SD1} will be enlarged as almost $\ast$-Ricci-Bourguignon Soliton\cite{DSP}.\\

On a Riemannian or Pseudo-Riemannian Manifold $(\mathfrak{M},g)$ on a vanishing vector field $\varrho$ is known as torse-forming \cite{KY}, whether
On the other hand,if 
\begin{equation}\label{1.5}
   \nabla_{\mathbb{Y}_{1}} \varrho=\Psi \mathbb{Y}_1 + \gamma(\mathbb{Y}_1) \varrho 
\end{equation}
where $\nabla$ is the Levi-Civita connection of $g$, $\Psi$ is a smooth function,and $\gamma$ is a $1$-form, then a non vanishing vector field $\varrho$ on a Riemannian or Pseudo-Riemannian manifold $(\mathfrak{M},g)$ is known to be torse-forming. Also in equation \eqref{1.5}, if the $1$-form disappears identically, the vector field is referred to as con-circular. If the function $\Psi=1$ and the $1$-form disappear identically in \eqref{1.5}, then the vector field is referred to as concurrent. If the function $\Psi=0$ in \eqref{1.5} holds true, then the vector field is considered as recurrent. Eventually, the vector field is referred to be a parallel vector field if in \eqref{1.5} $\Psi=\gamma=0$.\\

  Chen \cite{YBC} presented the torqued vector field, a new vector field , in 2017. A vector field is referred to be a torqued vector field if it satisfies equation \eqref{1.5}  with $\gamma(\varrho)=0$.\\
  \\
  The overview of this article will be accompanied as: In section $3$, we discuss about some results of Kenmotsu manifold equipped with Ricci Bourguignon Soliton. We confessed the nature of Ricci Bourguignon soliton with Kenmotsu manifold as compressing, balancing or enlarging accordingly. In addition, we have found some results with conformal killing vector field. In section $4$, we established few curvature properties of Kenmotsu manifold with Ricci Bourguignon Soliton. In section $5$, we find some results of Kenmotsu manifold with Ricci Bourguignon Soliton as torse-forming vector field. In section $5$, we derived an example a $5$-dimensional Kenmotsu manifold equipped with $\ast$-Ricci Bourguignon Soliton. 
\section{Preliminaries}\label{2}
\hspace{0.5cm}
We examine an almost contact metric manifold $\mathfrak{M}$ that is odd dimensional. It has an almost contact metric structure $(\Phi, \xi, \eta, g)$ where  $\xi$ is a vector field,$\Phi$ is a (1,1)- tensor field, $\eta$ is a 1-form and $g$ is the consistent Riemannian metric such that
\begin{equation}
\Phi^2(\mathbb{Y}_{1})=-\mathbb{Y}_{1}+\eta(\mathbb{Y}_{1})\xi, \eta\circ\Phi=0, \Phi\xi=0,\eta(\xi)=1,
\end{equation}
\begin{equation}
g(\Phi{\mathbb{Y}_{1}}, \Phi{\mathbb{Y}_{2}})=g(\mathbb{Y}_{1},\mathbb{Y}_{2})-\eta(\mathbb{Y}_{1})\eta(\mathbb{Y}_{2}),
\end{equation}
\begin{equation}
g(\mathbb{Y}_{1},\Phi{\mathbb{Y}_{2}})=-g(\Phi{\mathbb{Y}_{1}},\mathbb{Y}_{2}),
\end{equation}
\begin{equation}
g(\mathbb{Y}_{1},\xi)=\eta(\mathbb{Y}_{1}),
\end{equation}
for all the vector fields $\mathbb{Y}_{1}, \mathbb{Y}_{2} \in \chi{(\mathfrak{M})}$.\\
It is known that an almost contact metric manifold can be a Kenmotsu manifold if,
\begin{equation}
(\nabla_{\mathbb{Y}_{1}}\Phi)\mathbb{Y}_{2}=-g(\mathbb{Y}_{1},\Phi{\mathbb{Y}_{2}})\xi-\eta{(\mathbb{Y}_{2})}\Phi{\mathbb{Y}_{1}},
\end{equation}
\begin{equation}
\nabla_{\mathbb{Y}_{1}}\xi=\mathbb{Y}_{1}-\eta(\mathbb{Y}_{1})\xi,
\end{equation}
where $\nabla$ signifies the Riemannian connection of $g$.\\
The following relations are hold in a Kenmotsu manifold are,\\
\begin{equation}
    \eta(\mathbb{Y}_{1},\mathbb{Y}_{2})\mathbb{Y}_{3}=g(\mathbb{Y}_{1},\mathbb{Y}_{3})\eta(\mathbb{Y}_{2})-g(\mathbb{Y}_{2},\mathbb{Y}_{3})\eta(\mathbb{Y}_{1}),
\end{equation}
\begin{equation}
    R(\mathbb{Y}_{1},\mathbb{Y}_{2})\xi=\eta{(\mathbb{Y}_{1})}\mathbb{Y}_{2}-\eta{(\mathbb{Y}_{2})}\mathbb{Y}_{1},
\end{equation}
\begin{equation}
    R(\mathbb{Y}_{1},\xi)\mathbb{Y}_2=g(\mathbb{Y}_{1},\mathbb{Y}_{2})\xi-\eta{(\mathbb{Y}_{2})}\mathbb{Y}_{1},
\end{equation}
here $R$ stands for the Riemannian curvature tensor. \\
\begin{equation}
\mathbb{S}(\mathbb{Y}_{1},\xi)=-2n\eta(\mathbb{Y}_{1})
\end{equation}
\begin{equation}
    \mathbb{S}(\Phi{\mathbb{Y}_{1}},\Phi{\mathbb{Y}_{2}})=\mathbb{S}(\mathbb{Y}_{1},\mathbb{Y}_{2})+2n\eta(\mathbb{Y}_{1})\eta(\mathbb{Y}_{2}),
\end{equation}
\begin{equation}
   (\nabla_{\mathbb{Y}_{1}}\eta)\mathbb{Y}_{2}=g(\mathbb{Y}_{1},\mathbb{Y}_{2})-\eta{(\mathbb{Y}_{1})}\eta{(\mathbb{Y}_{2})}, 
\end{equation}
for every vector fields $\mathbb{Y}_{1}, \mathbb{Y}_{2}, \mathbb{Y}_{3} \in \chi(\mathfrak{M})$.\\
we know that,\\
\begin{equation}
    (\mathscr{L}_\xi{g})(\mathbb{Y}_{1},\mathbb{Y}_{2})=g(\nabla_{\mathbb{Y}_{1}}\xi,\mathbb{Y}_{2})+g(\mathbb{Y}_{1},\nabla_{\mathbb{Y}_{2}}\xi),
\end{equation}
Similarly, for every vector fields $\mathbb{Y}_{1},\mathbb{Y}_{2} \in \chi(\mathfrak{M})$.\\
Then, we using (2.6) and (2.13), we get \\
\begin{equation}
    (\mathscr{L}_\xi{g})(\mathbb{Y}_{1},\mathbb{Y}_{2})=2[g(\mathbb{Y}_{1},\mathbb{Y}_{2})-\eta(\mathbb{Y}_{1})\eta(\mathbb{Y}_{2})]
\end{equation}
\medskip
    The $\ast$-Ricci tensor on an odd - dimensional Kenmotsu manifold \cite{RA} is
    \begin{equation}\label{2.15}
        \mathbb{S}^{\ast}(\mathbb{Y}_{1},\mathbb{Y}_{2}) = \mathbb{S}(\mathbb{Y}_{1},\mathbb{Y}_{2}) + (2n-1) g(\mathbb{Y}_{1},\mathbb{Y}_{2}) + \eta(\mathbb{Y}_{1})\eta(\mathbb{Y}_{2}).
    \end{equation}
    In above equation, $\mathbb{Y}_{1}=\mathfrak{e}_i$ and $\mathbb{Y}_{2}=\mathfrak{e}_i$ represent a local orthonormal frame and we sum over $i$ varies from $1$ to $2n+1$, we obtain
    \begin{equation}\label{2.16}
        r^{\ast} = r + 4n^2,
    \end{equation}
    where $r^{\ast}$ is referred as the $\ast$- scalar curvature of $\mathfrak{M}$.

\section{Certain results concerning Kenmotsu Manifold equipped with $\ast$- RB Soliton}\label{3}
\hspace{0.5cm}
 Consider $\mathfrak{M}$,  an odd dimensional Kenmotsu Manifold where metric $g$ assures the $\ast$-Ricci Bourguignon Soliton $(g,\nu,\Omega,\omega)$. We utilize \eqref{2.15} and \eqref{2.16} in \eqref{1.4}, we obtain,
 \begin{equation}\label{3.1}
     (\mathscr{L}_ z g)(\mathbb{Y}_{1},\mathbb{Y}_{2})+2\mathbb{S}(\mathbb{Y}_{1},\mathbb{Y}_{2})\\=\\ 2[\Omega - (2 n-1)+4\omega n^2]g(\mathbb{Y}_{1},\mathbb{Y}_{2})\\ +2\omega r g(\mathbb{Y}_{1},\mathbb{Y}_{2})-2\eta(\mathbb{Y}_{1}) \eta(\mathbb{Y}_{2})
 \end{equation}
 for all vector fields $\mathbb{Y}_1,\mathbb{Y}_{2}$ on $\mathfrak{M}$. Then by reduction we obtain, 
 \begin{equation}\label{3.2}
 (\mathscr{L}_z g)(\mathbb{Y}_{1},\mathbb{Y}_{2})+2\mathbb{S}(\mathbb{Y}_{1},\mathbb{Y}_{2})=\\ 2\Lambda g(\mathbb{Y}_{1},\mathbb{Y}_{2})+2\omega r g(\mathbb{Y}_{1},\mathbb{Y}_{2})+2\mu \eta(\mathbb{Y}_{1}) \eta(\mathbb{Y}_{2})
 \end{equation}
Hence we can state
\begin{thm}\label{Th1}
   Suppose the metric $g$ of odd dimensional Kenmotsu Manifold assures $\ast$-Ricci Bourguignon soliton $(g,z,\Omega,\omega)$ then the soliton changes to an $\eta$-RB Soliton $(g,z,\Lambda,\mu,\omega)$ where 
   \begin{equation}\label{3.3}
   \Lambda=\Omega-(2n-1)+4\omega n^2
   \end{equation}
   \begin{equation}\label{3.4}
   \mu=-1.
   \end{equation}
\end{thm}
Let $z=\xi$,$\mathbb{Y}_{2}=\xi$ in $(3.2)$ and using $(2.14)$,$(2.4)$ and $(2.10)$, we get,
\begin{equation}\label{3.5}
\Omega= -\omega(4n^2+r)
\end{equation}
which gives
\begin{equation}\label{3.6}
r=-4n^2-\frac{\Omega}{\omega}
\end{equation}
where $\omega \neq 0$.

\begin{cor}
Let the metric $g$ of an odd dimensional Kenmotsu Manifold assures $\ast$-RB Soliton $(g,z,\Lambda,\Omega,\omega)$ then the soliton will be compressing, balancing or enlarging accordingly as, 
$\omega ( 4n^2+r)\lesseqqgtr0$ and this scalar curvature $r$ becomes $-4n^2-\frac{\Omega}{\omega}$ with $\omega \neq 0$.
\end{cor} 
Let us consider a $\ast$ - RB Soliton $(g,V,\Omega,\omega)$ on an odd dimensional Kenmotsu manifold $\mathfrak{M}$ as,
\begin{equation}\label{3.7}
\mathscr{L}_z g(\mathbb{Y}_{1},\mathbb{Y}_{2})+2\mathbb{S}^{\ast}(\mathbb{Y}_{1},\mathbb{Y}_{2})=2\Omega g(\mathbb{Y}_{1},\mathbb{Y}_{2})+2\omega r^{\ast}g(\mathbb{Y}_{1},\mathbb{Y}_{2})
\end{equation}
for all vector fields $\mathbb{Y}_{1},\mathbb{Y}_{2} \in \chi(\mathfrak{M})$.\\
Taking $\mathbb{Y}_{1}=\mathbb{Y}_{2}=\mathfrak{e}_i$ in $\eqref{3.7}$, here $\mathfrak{e}_i's$ are local orthonormal basis, and summate over $i=1,2,....,(2n+1)$ and using $(2.16)$, we obtain
\begin{equation}\label{3.8}
 div z-\Omega (2n+1)-[\omega(2n+1)-1][r+4n^2]=0 
\end{equation}
As the vector field $z$ is the gradient type, i.e., $z=grad(f)$, for $f$ is a smooth function on $\mathfrak{M}$, then the above equation $\eqref{3.8}$ will be reduced as,
\begin{equation}\label{3.9}
\bigtriangleup(f)=\Omega (2n+1)+ [\omega(2n+1)-1][r+4n^2]
\end{equation}
here $\bigtriangleup(f)$ is the Laplacian equation satisfied by $f$.
Hence we can state,
\begin{thm}
If the metric $g$ of an odd dimensional Kenmotsu manifold satisfies $\ast$-RB Soliton $(g,z,\Omega,\omega)$ where $z$ is the gradient of smooth function $f$ then the Laplacian equation satisfied by $f$ is \\
$\bigtriangleup f=\Omega (2n+1)+[\omega(2n+1)-1][r+4n^2] $
\end{thm}
Suppose the vector field is solenoidal, i.e.,$div z=0$ then (3.8) becomes 
\begin{equation}\label{3.10}
    r=\frac{\Omega(2n+1)}{[1-\omega(2n+1)]}-4n^2
\end{equation}
Inversely, assume the scalar curvature of the manifold
$ r=\frac{\Omega(2n+1)}{[1-\omega(2n+1)]}-4n^2$.
Then from (3.8), we have $div z=0$, which infers that $z$
is Solenoidal.
This completes the result.
Hence we can write,
\begin{rem}
 Let the metric $g$ and of an odd dimensional Kenmotsu manifold satisfy the $\ast$-RB soliton $(g,\nu,\Omega,\omega)$. Then  $z$(vector field) is
 Solenoidal if and only if the $r$(scalar curvature) will be
$\frac{\Omega(2n+1)}{[1-\omega(2n+1)]}-4n^2$.
\end{rem}
\begin{defn}
    A vector field $z$ is known to be a conformal killing vector field iff the following expression satisfies,
\begin{equation}\label{3.11}
    (\mathscr{L}_z{g})(\mathbb{Y}_{1},\mathbb{Y}_{2})=2\Omega{g}(\mathbb{Y}_{1},\mathbb{Y}_{2})
\end{equation}
where $\Omega$ is conformally scalar function of the coordinates.\\
\end{defn}
  Furthermore, the conformal killing vector field $z$ will be demonstrated as proper if the $ \Omega$ is not constant. If $\Omega$ is constant, then $z$ is known as homothetic vector field, If the constant $\Omega$ convert as non-zero, then $z$ is known as proper homothetic vector field. If $\Omega=0$ in the equation (3.5), then $z$ is said to be a killing vector field \cite{RA}.\\
  
  Using the definition (3.5), we can say, \\ 

     Let $(g, \omega, z, \Omega)$ be a $\ast$-RB soliton on an odd dimensional Kenmotsu manifold $\mathfrak{M}$, where $z$ is a conformal killing vector field. Then using (1.4), (2.16), (3.11) and setting $\mathbb{Y}_{2}=\xi$, we have
\begin{equation}\label{3.12}
    [\alpha-\omega(r+4n^2)-\Omega]\eta(\mathbb{Y}_1)=0
\end{equation}
      since $\eta(\mathbb{Y}_{1})\neq{0}$
\begin{equation}\label{3.13}
    \alpha=\omega(r+4n^2)+\Omega
\end{equation}
Hence we can write,2
\begin{thm}
    Suppose the metric $g$ of an odd dimensional Kenmotsu manifold assures the $\ast$-RB soliton $(g,\omega, z, \Omega)$ where $z$ represents a conformal killing vector field then $z$ is
    \begin{enumerate}
        \item  defined as proper if the expression $\omega(r+4n^2)+\Omega$ is non constant.
        \item defined as Homothetic vector field if the expression $\omega(r+4n^2)+\Omega$ is constant.
        \item defined as Proper homothetic vector field if the expression $\omega(r+4n^2)+\Omega$ is non-zero constant.
        \item  defined as Killing vector field if the expression is $-\omega(r+4n^2)$.
    \end{enumerate}
\end{thm}
\section{Some well-known curvature properties on Kenmotsu manifold satisfying $\ast$-RB soliton}\label{4}
G. P. Pokhariyal and R. S. Mishra  \cite{MP} discovered a new curvature tensor and they investigated its characteristics in 1970. This $(0,4)$ type
tensor in an $n$-dimensional Riemannian manifold, denoted as
$\mathcal{W}_2$, is defined by,

\begin{equation} \label{4.1}
    \mathcal{W}_2(\mathbb{Y}_{1},\mathbb{Y}_{2},\mathbb{Y}_{3},\mathbb{Y}_{4})=\grave{R}(\mathbb{Y}_{1},\mathbb{Y}_{2},\mathbb{Y}_{3},\mathbb{Y}_{4})+\frac{1}{n-1}[g(\mathbb{Y}_{1},\mathbb{Y}_{3})\mathbb{S}(\mathbb{Y}_{2},\mathbb{Y}_{4})-g(\mathbb{Y}_{2},\mathbb{Y}_{3})\mathbb{S}(\mathbb{Y}_{1},\mathbb{Y}_{4})],
\end{equation}

for all vector fields $\mathbb{Y}_{1}, \mathbb{Y}_{2},\mathbb{Y}_{3}, \mathbb{Y}_{4}$ and $\grave{\mathfrak{R}}(\mathbb{Y}_{1},\mathbb{Y}_{2},\mathbb{Y}_{3},\mathbb{Y}_{4})=g(\mathfrak{R}(\mathbb{Y}_{1},\mathbb{Y}_{2})\mathbb{Y}_{3},\mathbb{Y}_{4})$.
\begin{defn}
A Riemannian manifold is considered to be $\mathcal{W}_2$-flat if its $\mathcal{W}_2$-curvature tensor disappears identically.
\end{defn}

Let us assume that we have a $\mathcal{W}_2$-flat Kenmotsu manifold.\\
Then from equation \eqref{4.1} and the definition of $\mathcal{W}_2$-flat, we have,

\begin{equation} \label{4.2}
  \grave{\mathfrak{R}}(\mathbb{Y}_{1},\mathbb{Y}_{2},\mathbb{Y}_{3},\mathbb{Y}_{4})+\frac{1}{2n}[g(\mathbb{Y}_{1},\mathbb{Y}_{3})\mathbb{S}(\mathbb{Y}_{2},\mathbb{Y}_{4})-g(\mathbb{Y}_{2},\mathbb{Y}_{3})\mathbb{S}(\mathbb{Y}_{1},\mathbb{Y}_{4})]=0.
\end{equation}

We substitute $\mathbb{Y}_{1} = \mathbb{Y}_{4} = \mathfrak{e}_i$ in $\eqref{4.2}$, where $\mathfrak{e}_i$'s are a local orthonormal basis and $\sum_{i=1}^{2n+1}$ we get,

\begin{equation} \label{4.3}
  \mathbb{S}(\mathbb{Y}_{2},\mathbb{Y}_{3})=\frac{r}{2n+1}g(\mathbb{Y}_{2},\mathbb{Y}_{3}).
\end{equation}
Taking the help of \eqref{4.3}, \eqref{3.2} takes,
\begin{equation}\label{4.4}
 (\mathscr{L}_z g)(\mathbb{Y}_{1},\mathbb{Y}_{2})+\frac{2r}{2n+1}g(\mathbb{Y}_{1},\mathbb{Y}_{2})=\\ 2\Lambda g(\mathbb{Y}_{1},\mathbb{Y}_{2})+2\omega r g(\mathbb{Y}_{1},\mathbb{Y}_{2})+2\mu \eta(\mathbb{Y}_{1}) \eta(\mathbb{Y}_{2})
 \end{equation}
 Now if the vector field $z$ is concircular vector field \cite{CBY}, then we know for any smooth function $z$,

\begin{equation}\label{4.5}
  \nabla_{\mathbb{Y}_{1}} z=z \mathbb{Y}_{1},
\end{equation}
for all vector fields $\mathbb{Y}_{1}$.
Making use the rule of Lie derivative and \eqref{4.5}, one can have,

\begin{equation} \label{4.6}
  (\mathscr{L}_z g)(\mathbb{Y}_{2},\mathbb{Y}_{3})=2z g(\mathbb{Y}_{2},\mathbb{Y}_{3}).
\end{equation}
Thus using the preceding equation, \eqref{4.4} takes the form,
\begin{equation}\label{4.7}
    \Bigl[\Lambda+\omega r-\frac{r}{2n+1}-\nu \Bigl]g(\mathbb{Y}_{1},\mathbb{Y}_{2})=-\mu\eta(\mathbb{Y}_{1})\eta(\mathbb{Y}_{2})
\end{equation}
Taking $\mathbb{Y}_{2}=\xi$, we get,
\begin{equation}\label{4.8}
    \Lambda+\omega r-\frac{r}{2n+1}-\nu+\mu=0,
\end{equation}
as $\eta(\mathbb{Y}_{1}) \neq 0, \text{for any vector field}$ $\mathbb{Y}_{1}$.\\
From Theorem \ref{Th1}, we can substitute the values of $\Lambda$ and $\mu$, so we can shape \eqref{4.8} as,
$$\Omega=z+\frac{r}{2n+1}+2n-\omega r-4\omega n^2.$$
Hence we can have,
\begin{thm}
Suppose the metric $g$ of an odd dimensional $\mathcal{W}_2$-flat Kenmotsu manifold equipped with $\ast$-RB Soliton $(g,z,\Omega,\omega)$ where $z$ is concircular vector field, then the soliton will be compressing, balancing, enlarging accordingly  as, $z+\frac{r}{2n+1}+2n-\omega r-4\omega n^2\gtreqqless 0$
\end{thm}

In $n$-dimensional Riemannian manifold, Mantica \cite{YJCA} and Suh invented a new curvature tensor which is mentioned as $Q$ and it is stated as,
\begin{equation}\label{4.9}
Q(\mathbb{Y}_{1},\mathbb{Y}_{2})\mathbb{Y}_{3}=\mathfrak{R}(\mathbb{Y}_{1},\mathbb{Y}_{2})\mathbb{Y}_{3}-\frac{\psi}{n-1}[g(\mathbb{Y}_{2},\mathbb{Y}_{3})X-g(\mathbb{Y}_{1},\mathbb{Y}_{3})\mathbb{Y}_{2}]
\end{equation}
here $\psi$ is an arbitrary scalar function. Then tensor $Q$ is said to be a $Q$-curvature tensor.\\
\begin{defn}
A Riemannian manifold is called as $Q-flat$ \cite{SSS} if its $Q$-curvature tensor disappears equivalently.
\end{defn}
Let us assume that we have $Q-flat$ Kenmotsu manifold. \\Then from equation \eqref{4.9} and the definition of $Q-flat$, we have,
\begin{equation}\label{4.10}
\grave{\mathfrak{R}}(\mathbb{Y}_{1},\mathbb{Y}_{2},\mathbb{Y}_{3},\mathbb{Y}_{4})= \frac{\psi}{2n} [g(\mathbb{Y}_{2},\mathbb{Y}_{3})g(\mathbb{Y}_{1},\mathbb{Y}_{4})-g(\mathbb{Y}_{1},\mathbb{Y}_{3})g(\mathbb{Y}_{2},\mathbb{Y}_{4})]
\end{equation}
By taking $\mathbb{Y}_{1} =\mathbb{Y}_{4} =\mathfrak{e}_i$, in \eqref{4.10}, where $\mathfrak{e}_i's$ are orthonormal basis and we summate over $i=1,2,..(2n+1)$ we have,
\begin{equation}\label{4.11}
\mathbb{S}(\mathbb{Y}_{2},\mathbb{Y}_{3}) = \psi g(\mathbb{Y}_{2},\mathbb{Y}_{3})
\end{equation}
The above equation is known as Einstein.\\
Then substitute \eqref{4.11} in \eqref{3.2}, we get,
\begin{equation}\label{4.12}
 (\mathscr{L}_z g)(\mathbb{Y}_{1},\mathbb{Y}_{2}) + 2 \psi g(\mathbb{Y}_{1},\mathbb{Y}_{2})=2\Lambda g(\mathbb{Y}_{1},\mathbb{Y}_{2})+2 \omega rg(\mathbb{Y}_{1},\mathbb{Y}_{2})+2 \mu \eta(\mathbb{Y}_{1}) \eta(\mathbb{Y}_{2})
\end{equation}
Taking $ \mathbb{Y}_{1} = \mathbb{Y}_{2} = \mathfrak{e}_i $ and we summate over $i=1,2,...,(2n+1)$,then 
\begin{equation}
    divz + (2n+1)(\psi - \Lambda -\omega r)-\mu=0
\end{equation}
\begin{equation}
    divz= \mu-(2n+1)(\psi-\Lambda-\omega r)
\end{equation}
If $z$ is Solenoidal,
\begin{equation}
    \mu-(2n+1)(\psi-\Lambda-\omega r)=0
\end{equation}
\begin{equation}
    \Lambda=\psi-\frac{\mu}{2n+1}-\omega r
\end{equation}
Using (3.3) and (3.4), we get
\begin{equation}
    \Omega-(2n+1)+4\omega n^2=\psi-\frac{(-1)}{2n+1}-\omega r
\end{equation}
\begin{equation}
    \Omega=\psi+\frac{1}{2n+1}+(2n-1)-\omega(r+4n^2)
\end{equation}
Hence we can state
\begin{thm}
    Suppose the metric $g$ of an odd dimensional $Q$-flat Kenmotsu manifold assures $\ast$-Ricci Bourguignon soliton $(g, z, \Lambda, \Omega, \omega)$ then the soliton will be compressing, balancing or enlarging if $$\psi+\frac{1}{2n+1}+(2n+1)-\omega(r+4n^2)\gtreqqless 0$$ provided $z$ is  Solenoidal.
\end{thm}
\section{Results on torse forming vector field on Kenmotsu manifold equipped with $\ast$-Ricci Bourguignon soliton}\label{5}
  Here, we make use of \eqref{1.5}, and we write,
\begin{equation}\label{5.1}
    (\mathcal{L}_\varrho g)(\mathbb{Y}_{1},\mathbb{Y}_{2})=g(\nabla_{\mathbb{Y}_{1}}  \varrho,\mathbb{Y}_{2})+g(\mathbb{Y}_{1},\nabla_{\mathbb{Y}_{2}} \varrho)
\end{equation}
Then we have,
\begin{equation}\label{5.2}
    (\mathscr{L}_\varrho g)(\mathbb{Y}_{1},\mathbb{Y}_{2})=2 \Psi g(\mathbb{Y}_{1},\mathbb{Y}_{2})+\gamma(\mathbb{Y}_{1}) g(\varrho,\mathbb{Y}_{2})+\gamma(\mathbb{Y}_{2}) g(\varrho,\mathbb{Y}_{1})
\end{equation}
We take $v=\varrho$ in \eqref{3.1} and using \eqref{5.2}, we obtain,
\begin{multline}\label{5.3}
    2[- \Psi + \Omega - (2n-1) +4\omega n^2 +\omega r]g(\mathbb{Y}_{1},\mathbb{Y}_{2})-2\mathbb{S}(\mathbb{Y}_{1},\mathbb{Y}_{2})-2\eta(\mathbb{Y}_{1})\eta(\mathbb{Y}_{2}) \\
     =\gamma(\mathbb{Y}_{1})g(\varrho,\mathbb{Y}_{2})+\gamma(\mathbb{Y}_{2})g(\varrho,\mathbb{Y}_{1})
\end{multline}
Now taking $\mathbb{Y}_{1},\mathbb{Y}_{2}=\mathfrak{e}_i$ and summing over $1$ to $(2n+1)$, we have,
\begin{equation}\nonumber
    [\Omega - \Psi - (2n-1)+\omega(r+4n^2)](2n+1)-r-1=r(\varrho)
\end{equation}
At the end, we obtain,
\begin{equation}\label{5.4}
 \Omega = \frac{r(\varrho)+r+1}{2n+1} + \Psi + (2n-1) -\omega(r+4n^2)
\end{equation}
Hence we can state,
\begin{thm}
    Suppose the metric $g$ of an odd dimensional Kenmotsu manifold acquires a $\ast$-Ricci Bourguignon soliton, where $\varrho$ is a torse-forming vector field, then
    $\Omega =\frac{r(\varrho)+r+1}{2n+1} + \Psi +(2n-1) - \omega(r+4n^2)$ the soliton is said to be compressing, balancing or enlarging accordingly as, $\frac{r(\varrho)+r+1}{2n+1} + \Psi +(2n-1) - \omega(r+4n^2)$ $\gtreqqless$ $0$ respectively.
\end{thm}
If the $\omega$-$1$-form vanishes identically in \eqref{5.4}, then 
    $\Omega = \frac{r(\varrho)+r+1}{2n+1} + \Psi + (2n-1)$.\\
      In \eqref{5.4}, if $\omega$ $1$-form vanishes identically and the function $\Psi=1$ then, $\Omega = \frac{r(\varrho)+r+1}{2n+1} + (2n-1)$.\\
      If $\Psi = 0$, then \eqref{5.4} becomes $\Omega = \frac{r(\varrho)+r+1}{2n+1} + (2n-1) - \omega(r+4n^2)$. \\
     If $\gamma(\varrho)=0$, then \eqref{5.4} becomes, $\Omega = \frac{r+1}{2n+1} + \Psi +(2n-1) - \omega(r+4n^2)$.\\
 At the end we have, if $\Psi,\omega = 0$, then \eqref{5.4}
will be, $\Omega=\frac{r(\varrho)+r+1}{2n+1}+(2n-1)$.
It follows the Corollary,
\begin{cor}
    Suppose the metric $g$ of a odd dimensional Kenmotsu manifold assures $\ast$-Ricci Bourguignon soliton, where $\varrho$ is a torse-forming vector field, then
    \begin{enumerate}
    \item It is $concircular$, $\Omega = \frac{r(\varrho)+r+1}{2n+1} + \Psi + (2n-1)$ and the soliton is said to be compressing, balancing or enlarging with respect to $\frac{r(\varrho)+r+1}{2n+1} + \Psi + (2n-1)$ $\gtreqqless 0$. \\
    \item It is $concurrent$, $\Omega = \frac{r(\varrho)+r+1}{2n+1} + (2n-1)$ and the soliton is known to be compressing, balancing or enlarging with respect to $\frac{r(\varrho)+r+1}{2n+1} + (2n-1) $ $\gtreqqless 0$ respectively.\\
    \item  It is $recurrent$, $\Omega = \frac{r(\varrho)+r+1}{2n+1} + (2n-1) - \omega(r+4n^2)$ and the soliton is said to be compressing, balancing or enlarging accordingly as $\frac{r(\varrho)+r+1}{2n+1} + (2n-1) - \omega(r+4n^2)$ $\gtreqqless 0$.\\
    \item It is $torqued$, $\Omega = \frac{r+1}{2n+1} + \Psi +(2n-1) - \omega(r+4n^2)$ and the soliton is said to be compressing, balancing or enlarging with respect to $ \frac{r+1}{2n+1} + \Psi +(2n-1) - \omega(r+4n^2)$ $\gtreqqless 0$. \\
    \item It is $parallel$, $\Omega=\frac{r(\varrho)+r+1}{2n+1}+(2n-1)$ and the soliton is known to be compressing, balancing or enlarging accordingly as $\frac{r(\varrho)+r+1}{2n+1}+(2n-1)$ $\gtreqqless 0$ respectively.
    \end{enumerate}
\end{cor}
\section{Examples of a $5$-dimensional Kenmotsu manifold admitting $\ast$-Ricci Bourguignon soliton}\label{6}

\hspace{0.5cm} Let us assume that the manifold $\mathfrak{M}$ to be $5$ dimensional. \\
Here $\mathfrak{M}={(\mathbb{Y}_{1},\mathbb{Y}_{2},\mathbb{Y}_{3},\mathbb{Y}_{4},\mathbb{Y}_{5}) \in \mathbb{R}^5}$ \cite{RID,SS1} where $(\mathbb{Y}_{1},\mathbb{Y}_{2},\mathbb{Y}_{3},\mathbb{Y}_{4},\mathbb{Y}_{5})$ are the coordinates in $\mathbb{R}^5$. The vector field defined here are: \\
\[
\mathfrak{e}_1 = \mathfrak{e}^{-v} \frac{\partial}{\partial \mathbb{Y}_{1}}, \quad
\mathfrak{e}_2 = \mathfrak{e}^{-v} \frac{\partial}{\partial \mathbb{Y}_{2}}, \quad
\mathfrak{e}_3 = \mathfrak{e}^{-v} \frac{\partial}{\partial \mathbb{Y}_{3}}, \quad
\mathfrak{e}_4 = \mathfrak{e}^{-v} \frac{\partial}{\partial \mathbb{Y}_{4}}, \quad
\mathfrak{e}_5 = \frac{\partial}{\partial \mathbb{Y}_{5}}
\]
are linearly independent at each point of $\mathfrak{M}$. we describe  the metric $g$ as, \\
\begin{equation}\nonumber
 g(\mathfrak{e}_i,\mathfrak{e}_j) = 
\begin{cases} 
1, & \text{if}  \ i=j \ \text{and} \ i,j \in (1,2,3,4,5),\\
0, & \text{otherwise.}
\end{cases}
\end{equation}
Consider $\eta$ be a $1-$form denoted by $\eta(\mathbb{Y}_{1})=g(\mathbb{Y}_{1},\mathfrak{e}_5)$ for any $\mathbb{Y}_{1} \in \chi(\mathfrak{M})$. we set $(1,1)$-tensor field $\Phi$ as:\\
\[
\Phi(\mathfrak{e}_5) = 0, \quad
\Phi(\mathfrak{e}_4) = -\mathfrak{e}_2, \quad
\Phi(\mathfrak{e}_3) = -\mathfrak{e}_1, \quad
\Phi(\mathfrak{e}_2) = \mathfrak{e}_4, \quad
\Phi(\mathfrak{e}_1) =\mathfrak{e}_3 .
\]

Then it assure the relations $\eta(\xi)=1$, $ \Phi ^2(\mathbb{Y}_{1})=- \mathbb{Y}_{1}+ \eta(\mathbb{Y}_{1}) \xi$ and $g(\Phi \mathbb{Y}_1,\Phi \mathbb{Y}_2)=g(\mathbb{Y}_{1},\mathbb{Y}_{2})-\eta(\mathbb{Y}_{1})\eta(\mathbb{Y}_{2})$, where $\xi=\mathfrak{e}_5$ and $\mathbb{Y}_{1},\mathbb{Y}_{2}$ is any vector field on $\mathfrak{M}$. Therefore $(\mathfrak{M},\Phi,\xi,\eta,g)$ describes an almost contact structure on $\mathfrak{M}$.\\
We can now deduce that:
\begin{align*}
[\mathfrak{e}_5, \mathfrak{e}_1] &= -\mathfrak{e}_1, & [\mathfrak{e}_5, \mathfrak{e}_2] &= -\mathfrak{e}_2, & [\mathfrak{e}_5, \mathfrak{e}_3] &= -\mathfrak{e}_3, & [\mathfrak{e}_5, \mathfrak{e}_4] &=- \mathfrak{e}_4, \\
[\mathfrak{e}_4, \mathfrak{e}_1] &= 0, & [\mathfrak{e}_4, \mathfrak{e}_2] &= 0, & [\mathfrak{e}_4, \mathfrak{e}_3] &= 0, & [\mathfrak{e}_4, \mathfrak{e}_5] &= \mathfrak{e}_4, \\
[\mathfrak{e}_3, \mathfrak{e}_1] &= 0, & [\mathfrak{e}_3, \mathfrak{e}_2] &= 0, & [\mathfrak{e}_3, \mathfrak{e}_4] &= 0, & [\mathfrak{e}_3, \mathfrak{e}_5] &= \mathfrak{e}_3, \\
[\mathfrak{e}_2, \mathfrak{e}_1] &= 0, & [\mathfrak{e}_2, \mathfrak{e}_3] &= 0, & [\mathfrak{e}_2, \mathfrak{e}_4] &= 0, & [\mathfrak{e}_2, \mathfrak{e}_5] &= \mathfrak{e}_2, \\
[\mathfrak{e}_1, \mathfrak{e}_2] &= 0, & [\mathfrak{e}_1, \mathfrak{e}_3] &= 0, & [\mathfrak{e}_1, \mathfrak{e}_4] &= 0, & [\mathfrak{e}_1, \mathfrak{e}_5] &= \mathfrak{e}_1
\end{align*}
Consider $\nabla$ as the Levi-Civita connection associated with $g$.  According to Koszul's formula for any  $X, Y, Z \in \chi(\mathfrak{M}) $ is written as: \\
\[
2g(\nabla_{\mathbb{Y}_{1} \mathbb{Y}_{2}, \mathbb{Y}_{3}}) = \mathbb{Y}_{1}g(\mathbb{Y}_{2}, \mathbb{Y}_{3}) + \mathbb{Y}_{2}g(\mathbb{Y}_{3}, \mathbb{Y}_{1}) - \mathbb{Y}_{3}g(\mathbb{Y}_{1}, \mathbb{Y}_{2}) - g(\mathbb{Y}_{1}, [\mathbb{Y}_{2}, \mathbb{Y}_{3}]) - g(\mathbb{Y}_{2}, [\mathbb{Y}_{1}, \mathbb{Y}_{3}]) + g(\mathbb{Y}_{3}, [\mathbb{Y}_{1}, \mathbb{Y}_{2}]),
\]
we can have:
\[
\begin{aligned}
\nabla_{\mathfrak{e}_5} \mathfrak{e}_1 &=0, & \nabla_{\mathfrak{e}_5} \mathfrak{e}_2 &= 0, & \nabla_{\mathfrak{e}_5} \mathfrak{e}_3 &= 0, & \nabla_{\mathfrak{e}_5} \mathfrak{e}_4 &= 0, & \nabla_{\mathfrak{e}_5} \mathfrak{e}_5 &= 0, \\
\nabla_{\mathfrak{e}_4} \mathfrak{e}_1 &= 0, & \nabla_{\mathfrak{e}_4} \mathfrak{e}_2 &= 0, & \nabla_{\mathfrak{e}_4} \mathfrak{e}_3 &= 0, & \nabla_{\mathfrak{e}_4} \mathfrak{e}_4 &= -\mathfrak{e}_5, & \nabla_{\mathfrak{e}_4} \mathfrak{e}_5 &= \mathfrak{e}_4, \\
\nabla_{\mathfrak{e}_3} \mathfrak{e}_1 &= 0, & \nabla_{\mathfrak{e}_3} \mathfrak{e}_2 &= 0, & \nabla_{\mathfrak{e}_3} \mathfrak{e}_3 &= -\mathfrak{e}_5, & \nabla_{\mathfrak{e}_3} \mathfrak{e}_4 &= 0, & \nabla_{\mathfrak{e}_3} \mathfrak{e}_5 &= \mathfrak{e}_3, \\
\nabla_{\mathfrak{e}_2} \mathfrak{e}_1 &= 0, & \nabla_{\mathfrak{e}_2} \mathfrak{e}_2 &= -\mathfrak{e}_5, & \nabla_{\mathfrak{e}_2} \mathfrak{e}_3 &= 0, & \nabla_{\mathfrak{e}_2} \mathfrak{e}_4 &= 0, & \nabla_{\mathfrak{e}_2} \mathfrak{e}_5 &= \mathfrak{e}_2, \\
\nabla_{\mathfrak{e}_1} \mathfrak{e}_1 &= -\mathfrak{e}_5, & \nabla_{\mathfrak{e}_1} \mathfrak{e}_2 &= 0, & \nabla_{\mathfrak{e}_1} \mathfrak{e}_3 &= 0, & \nabla_{\mathfrak{e}_1} \mathfrak{e}_4 &= 0, & \nabla_{\mathfrak{e}_1} \mathfrak{e}_5 &= \mathfrak{e}_1.
\end{aligned}
\]

Hence,  $ (\nabla _{\mathbb{Y}_{1}} \Phi) \mathbb{Y}_{2}= g (\Phi \mathbb{Y}_{1}, \mathbb{Y}_{2}) \xi - \eta(\mathbb{Y}_{2}) \Phi(\mathbb{Y}_{1}) $ is reassured for arbitrary $\mathbb{Y}_{1},\mathbb{Y}_{2} \in \chi(\mathfrak{M})$. So $(\mathfrak{M},\Phi,\xi,\eta,g)$ will be a Kenmotsu manifold. \\
The non-zero curvature tensor elements are;
\begin{align*}
\mathfrak{R}(\mathfrak{e}_5, \mathfrak{e}_3) \mathfrak{e}_5 &= \mathfrak{e}_3, & \mathfrak{R}(\mathfrak{e}_5, \mathfrak{e}_4) \mathfrak{e}_5 &= \mathfrak{e}_4, & \mathfrak{R}(\mathfrak{e}_4, \mathfrak{e}_5)\mathfrak{e}_4 &= \mathfrak{e}_5, \\
\mathfrak{R}(\mathfrak{e}_3, \mathfrak{e}_4)\mathfrak{e}_4 &= -\mathfrak{e}_3, & \mathfrak{R}(\mathfrak{e}_3, \mathfrak{e}_5)\mathfrak{e}_3 &= \mathfrak{e}_5, & \mathfrak{R}(\mathfrak{e}_3, \mathfrak{e}_4 )\mathfrak{e}_3 &= \mathfrak{e}_4, \\
\mathfrak{R}(\mathfrak{e}_2, \mathfrak{e}_5)\mathfrak{e}_5 &= -\mathfrak{e}_2, & \mathfrak{R}(\mathfrak{e}_2, \mathfrak{e}_4)\mathfrak{e}_4 &=- \mathfrak{e}_2, & \mathfrak{R}(\mathfrak{e}_2, \mathfrak{e}_3)\mathfrak{e}_3 &= -\mathfrak{e}_2, \\
\mathfrak{R}(\mathfrak{e}_2, \mathfrak{e}_5)\mathfrak{e}_2 &= \mathfrak{e}_5, & R(\mathfrak{e}_2, \mathfrak{e}_4)\mathfrak{e}_2 &= \mathfrak{e}_4, & \mathfrak{R}(\mathfrak{e}_2, \mathfrak{e}_3)\mathfrak{e}_2 &= -\mathfrak{e}_3, \\
\mathfrak{R}(\mathfrak{e}_1, \mathfrak{e}_5)\mathfrak{e}_1 &= -\mathfrak{e}_5, & R(\mathfrak{e}_1, \mathfrak{e}_4)\mathfrak{e}_1 &= -\mathfrak{e}_4, & \mathfrak{R}(\mathfrak{e}_1, \mathfrak{e}_3)\mathfrak{e}_1 &= \mathfrak{e}_3, \\
\mathfrak{R}(\mathfrak{e}_1, \mathfrak{e}_2)\mathfrak{e}_1 &= \mathfrak{e}_2, & \mathfrak{R}(\mathfrak{e}_1, \mathfrak{e}_5)\mathfrak{e}_5 &= -\mathfrak{e}_1, & \mathfrak{R}(\mathfrak{e}_1, \mathfrak{e}_4)\mathfrak{e}_4 &= -\mathfrak{e}_1, \\
    \mathfrak{R}(\mathfrak{e}_1, \mathfrak{e}_3)\mathfrak{e}_3 &= -\mathfrak{e}_1, & \mathfrak{R}(\mathfrak{e}_1, \mathfrak{e}_2)\mathfrak{e}_2 &= - \mathfrak{e}_1.
\end{align*}
By the preceding results, we hold:
$
S(\mathfrak{e}_i, \mathfrak{e}_i) = -4 \quad \text{for } i = 1, 2, 3, 4, 5
$
and
\begin{equation}\label{6.1}
\mathbb{S}(\mathbb{Y}_{1}, \mathbb{Y}_{2}) = -4g(\mathbb{Y}_{1}, \mathbb{Y}_{2}) \quad \forall \mathbb{Y}_{1}, \mathbb{Y}_{2} \in \chi(\mathfrak{M}).
\end{equation}
Reducing the above, we have:
$
r = \sum_{i=1}^{5} \mathbb{S}(\mathfrak{e}_i, \mathfrak{e}_i) = -20 = -2n(2n + 1)
$
where the dimension of the manifold $(2n + 1 = 5)$.

In addition, we have:
$
\mathbb{S}^*(\mathfrak{e}_i, \mathfrak{e}_i) =
\begin{cases}
-1, & \text{if } i = 1, 2, 3, 4, \\
0, & \text{if } i = 5.
\end{cases}
$
Therefore, \\
\begin{equation}\label{6.2}
\mathbb{S}^{\ast}(\mathbb{Y}_{1}, \mathbb{Y}_{2}) = -g(\mathbb{Y}_{1}, \mathbb{Y}_{2}) + \eta(\mathbb{Y}_{1})\eta(\mathbb{Y}_{2}) \quad \forall \mathbb{Y}_{1}, \mathbb{Y}_{2} \in \chi(\mathfrak{M}).
\end{equation}
Hence,
\begin{equation}\label{6.3}
r^* = \text{Tr}(\mathbb{S}^{\ast}) = -4. \\
\end{equation}
Consider a vector field $z$ as \\
\begin{equation}\label{6.4}
z = \mathbb{Y}_{1}\frac{\partial}{\partial \mathbb{Y}_{1}} + \mathbb{Y}_{2}\frac{\partial}{\partial \mathbb{Y}_{2}} + \mathbb{Y}_{3}\frac{\partial}{\partial \mathbb{Y}_{3}} + \mathbb{Y}_{4}\frac{\partial}{\partial \mathbb{Y}_{4}} + \frac{\partial}{\partial \mathbb{Y}_{5}}
\end{equation}
Also, By the above results, we can justify that \\
\begin{equation}\label{6.5}
(\mathscr{L}_z g)(\mathbb{Y}_{1}, \mathbb{Y}_{2}) = 4\{g(\mathbb{Y}_{1}, \mathbb{Y}_{2}) - \eta(\mathbb{Y}_{1})\eta(\mathbb{Y}_{2})\}
\end{equation}
which adheres for all $\mathbb{Y}_{1}, \mathbb{Y}_{2} \in \chi(\mathfrak{M})$. Hence, considering \eqref{6.5}, we conclude \\
\begin{equation}\label{6.6}
\sum_{i=1}^{5} (\mathscr{L}_z g)(\mathfrak{e}_i, \mathfrak{e}_i)  =  16.
\end{equation}
Here, we take $\mathbb{F}_{1},\mathbb{F}_{2}=\mathfrak{e}_i$ in \eqref{1.4} and summing over $i$ varies from $1$ to $5$ and taking \eqref{6.3} and \eqref{6.6},
we obtain,
\begin{equation}\label{6.7}
\Omega = 4 \omega + \frac{4}{5}
\end{equation}
This $\Omega$ satisfies \eqref{3.8}. Therefore, $g$ satisfies $\ast$-Ricci Bourguignon soliton on Kenmotsu manifold $\mathfrak{M}$ which is of dimension $5$.

\section{Real world Implementations}
The exploration of $\ast$-Ricci Bourguignon on Kenmotsu manifold connects pure mathematics with applied theories like string theory, optical communication and general relativity. Solitons are generally known as solitary waves which has a vast applications in general relativity and cosmology.
\begin{itemize}
    \item Generalized Ricci solitons in Cosmological models: \\ Ricci soliton structure are self similar solution in the dynamics of spacetime. They offer realization to cosmic evolution and behavior of singularities. The study of contact structure offers a framework for learning such occurrences in an odd dimensional spacetimes.
    \item Existence of contact geometry in Classical Mechanics:\\ As a part of contact manifold, it has a origin in the mathematical representation of Classical Mechanics.
    \item Machine Learning and Contact Geometry:\\
    Kenmotsu Manifolds are used to describe high dimensional data manifolds in maching learning because of their contact geometry and soliton structures.
\end{itemize}

\textbf{Declaration of Competing Interest:}
The authors declare that they have no known competing financial or personal relationships that could have appeared to influence the work reported in this paper.


\end{document}